\documentclass[12pt]{amsart}
\usepackage{amsfonts}
\usepackage{latexsym}
\usepackage{amssymb}
\usepackage{amsmath}

\newcommand{\G}{\mathbb{G}}

\newcommand{\N}{\mathbb{N}}

\newcommand{\R}{\mathbb{R}}

\newcommand{\cG}{\mathcal{G}}
\newcommand{\cH}{\mathcal{H}}

\newcommand{\cL}{\mathcal{L}}

\newcommand{\cS}{\mathcal{S}}

\newcommand{\cp}{\mathfrak{p}}

\renewcommand{\span}{\mbox{span}}

\newcommand{\ep}{\varepsilon}

\newcommand{\sm}{\setminus}

\newcommand{\res}{\mbox{\LARGE{$\llcorner$}}}

\newcommand{\lra}{\longrightarrow}

\newcommand{\der}{\partial}

\newcommand{\Bx}{\mbox{\rm\small Box}}

\newtheorem{The}{Theorem}[section]
\newtheorem{Lem}{Lemma}

\newtheorem{Def}{Definition}
\newtheorem{Pro}{Proposition}[section]
\newtheorem{Cor}{Corollary}[section]

\begin{document}

\title
[Measure of curves in graded groups]
{{\bf Measure of curves in graded groups}}
\author{Riikka Korte}
\address{Riikka Korte, Department of Mathematics and Statistics \\ 
P.O. Box 68, FI-00014 University of Helsinki, Finland. }
\email{Riikka.Korte@helsinki.fi}
\author{Valentino Magnani}
\address{Valentino Magnani, Dipartimento di Matematica \\
Largo Bruno Pontecorvo 5 \\ I-56127, Pisa, Italy.}
\email{magnani@dm.unipi.it}

\begin{abstract}
We establish an area-type formula for the intrinsic spherical Hausdorff measure of
every regular curve embedded in an arbitrary graded group.
\end{abstract}

\maketitle

\tableofcontents

\footnoterule{
The second author has been supported by ''ERC ADG Grant GeMeTneES'' 
}

\pagebreak

\section{Introduction}

In the seminal paper by M. Gromov, \cite{Gr1}, an interesting formula for the Hausdorff 
dimension of submanifolds in arbitrary Carnot-Carath\'eodory spaces can be found, 
see Section 0.6 B. Once we have the Hausdorff dimension, the subsequent question
is to investigate whether the corresponding Hausdorff measure restricted to the submanifold
is positive, locally finite and can be computed by an integral formula. 

This is exactly the case, when a suitable ``negligibility condition'' is assumed
to be satisfied by the single submanifold, according to results of \cite{MagVit}.
In fact, here an area-type formula for $\cS^q\res\Sigma$ is established, where 
$q$ is the degree of the submanifold $\Sigma$ contained in a graded group $\G$ 
and $\cS^q$ represents the spherical Hausdorff measure
with respect to the fixed homogeneous distance of the group.
Notice that 0.6 B of \cite{Gr1} joined with Remark 4.2 of \cite{Mag12B} shows that
the degree of $\Sigma$ coincides with its Hausdorff dimension. In short, 
the validity of an area type formula for $\cS^q(\Sigma)$ is a consequence of the 
validity of the negligibility condition, that we recall in the following
\begin{Def}[Negligibility condition]\rm
Let $\Sigma$ be a submanifold of degree $q$, that is embedded in a stratified group $\G$.
Then $\Sigma$ satisfies the {\em negligibility condition} if 
\begin{equation}\label{negcon}
\cH^q\left(\left\{x\in\Sigma: d_\Sigma(x)<q\right\}\right)=0.
\end{equation}
\end{Def}
Here $d_\Sigma(x)$ is the pointwise degree of $\Sigma$ at $x$, \cite{MagVit}.
The validity of \eqref{negcon} for arbitrary submanifolds in arbitrary stratified groups 
is still an intriguing open question, where the regularity of the submanifold is an 
important parameter to be fixed and that depends on both the dimension of the submanifold
and on the step of the group. On the other hand, condition \eqref{negcon} holds
for $C^1$ smooth non-horizontal submanifolds, \cite{Mag6, Mag12A}, for $C^{1,1}$ 
smooth submanifolds of two step groups, \cite{Mag12B}, and for $C^{1,1}$ smooth 
submanifolds of the Engel group, \cite{LeDMag}. 
The proofs in these different cases rely on different approaches.

The present work continues this program, using the recent approach applied in
\cite{LeDMag}. Then we show that \eqref{negcon} also holds
for all $C^1$ smooth curves of an arbitrary graded group, namely, we have the following
\begin{The}[Negligibility]\label{negconcur}
Every 1-dimensional $C^1$ smooth submanifold contained in a graded group
satisfies the negligibility condition.
\end{The}
As it is clear from \cite{LeDMag}, proving negligibility requires 
some ad hoc arguments to treat the blow-up of different types of ``singular points'', 
that can occur in the study of the single submanifold.
This singular behaviour is somehow given by the interplay between the degree
of the submanifold $d(\Sigma)$ and the pointwise degree at the singular point $d_\Sigma(x)$.
In particular, this implies a precise differential contraint on the local
Taylor expansion of the submanifold near the singular point. 
In other terms, we show that all possible singularities of curves can be handled 
through a single argument by induction, that only requires $C^1$ regularity.

Notice that regularity is important, since the precision needed in the local description
of the submanifold depends on its ``local transversality''. 
Loosely speaking, in submanifolds of higher dimension low regularity 
allows for a larger size of the set of points with low degree. This fact has been
carefully studied in Heisenberg groups by Z. M. Balogh, \cite{Bal}.

As mentioned above, the negligibility condition implies an integral formula for the
spherical Hausdorff measure of the submanifold. This in turn follows by an intrinsic blow-up 
at each point of maximum degree, \cite{MagVit}. 
However, this blow-up requires $C^{1,1}$ smoothness, so that Theorem~\ref{negconcur} 
would only lead us to an area-type formula for 1-dimensional submanifolds of class $C^{1,1}$. 
This stronger regularity somehow conflicts with the $C^1$ smoothness
needed for the negligibility condition. 

In fact, the second issue of this paper concerns the possibility to perform an intrinsic blow-up
of curves at points of maximum degree, under $C^1$ regularity. 
The main point of the approach in \cite{MagVit} is that of finding a local 
reparametrization of the submanifold, using a family of curves on right neighbourhoods
of the origin that solve a Lipschitz ODEs on the submanifold.
For $C^1$ smooth submanifolds this approach fails, since the solutions of the same
continuous ODEs are no longer unique. Clearly, in the case of 1-dimensional submanifolds
two of these solutions are enough to parametrize a neighbourhood of the point of
maximum degree. This suggests to find a suitable blow-up reparametrization of the
curve around its point of maximum degree. This yields the correct blow-up estimates,
according to Proposition~\ref{blw}, and it allows us to establish the following 
\begin{The}\label{blwmea}
Let $\tilde g$ be any Riemannian metric on $\G$ and let $\Sigma$ be a $C^1$ 
smooth $1$-dimensional submanifold such that $d(\Sigma)=q$. Let $x\in\Sigma$ be 
such that $d_\Sigma(x)=q$. Then we have
\begin{equation}\label{1blw}
\lim_{r\to0^+}
\frac{\tilde\mu_1(\Sigma\cap B_{x,r})}{r^q}=
\frac{\theta(\tau_{\Sigma}^q(x))}{|\tau_{\Sigma}^q(x)|}\,,
\end{equation}
where $\tilde\mu_1$ is the Riemannian measure on $\Sigma$ with respect to the metric $\tilde g$.
\end{The}
Joining Theorem~\ref{negconcur} with Theorem~\ref{blwmea}, by standard differentiability 
theorems on measures, see \cite{Fed}, we get the following
\begin{Cor}
Let $\tilde g$ be a fixed Riemannian metric on $\G$. Let $\Sigma$ be a $C^1$ 
smooth one dimensional submanifold such that $d(\Sigma)=q$. Then we have
\begin{equation}\label{intcurv}
\int_{\Sigma}\theta(\tau^q_{\Sigma}(x))d\mathcal S^q(x)=\int_{\Sigma}|\tau^q_{\Sigma}(x)|\,d\tilde\mu_{1}(x).
\end{equation}
\end{Cor}
The possible anisotropy of the homogeneous distance $d$ makes the metric factor
$\theta(\tau^q_{\Sigma}(x))$ a function that depends on the direction $\tau^q_\Sigma(x)$, 
see Definition~\ref{metricf}.
Recall that, according to Subsection~\ref{subgv}, the stratified group $\G$ can be decomposed
into the direct sum of linear subspaces $H_j$ of degree $j$, where $j=1,\ldots,\iota$. 
Using this notation, we say that the distance $d$ has a {\em symmetry of degree $q$} if 
for all $\tau\in H_q$ we have 
\begin{eqnarray}
c_q=\theta(\tau)=\cH_{|\cdot|}^1\left(\span\{\tau\}\cap B_1\right),
\end{eqnarray}
then under this assumption the integral formula \eqref{intcurv} becomes
\begin{equation}\label{cntint}
c_q\,\cS^q(\Sigma)=\int_{\Sigma}|\tau^d_{\Sigma}(x)|\,d\tilde\mu_{1}(x).
\end{equation}
In connection with Gromov's dimension comparison probelm studied in \cite{BalTysWar09}, 
we wish to regard our work from a different perspective. Let $\G$ be a stratified group of step $\iota$. 
Then for any integer $q=1,\ldots,\iota$, one 
easily observes that the subgroup $L_q$ spanned by an element $e_i$ of degree $d_i=q$
has Hausdorff dimension equal to $q$. A consequence of our results is that these integers are 
the only possible Hausdorff dimensions of $C^1$ smooth curves.
In addition, each of these curves has an area-type formula \eqref{intcurv} 
that computes its spherical Hausdorff measure.

\section{Definitions and standard facts}\label{defstares}

A graded group $\mathbb G$ with topological dimension $n$ is a simply connected nilpotent 
Lie group with Lie algebra $\mathcal G$ having the grading $\mathcal G=V_{1}\oplus\ldots\oplus V_\iota$,
that satisfies the conditions $[V_i,V_j]\subset V_{i+j}$ for all $i,j\geq 1$, where 
$V_j=\{0\}$ whenever $j>\iota$. The integer $\iota$ is called the {\em step} of $\mathbb G$. 
Recall that the family of graded groups strictly contains the well known family of stratified groups.
The grading allows for defining dilations $\delta_r:\cG\rightarrow\cG$ as follows
\[
\delta_r\Big(\sum_{j=1}^\iota v_{j}\Big)=\sum_{j=1}^\iota r^jv_{j},\quad \mbox{for all} \quad r>0,
\]
where we have set $v=\sum_{j=1}^\iota v_{j}$ and $v_{j}\in V_{j}$.
Note that simply connected nilpotent Lie groups are diffeomorphic to their Lie algebra through the exponential mapping $\exp:\mathcal G\rightarrow \mathbb G$, hence dilations are automatically defined as group 
isomorphisms of $\mathbb G$ and will be denoted by the same symbol $\delta_r$.

We say that $\rho$ is a homogeneous distance on $\mathbb G$ if it is a continuous distance of $\mathbb G$ satisfying the following conditions:
\[
d(zx,zy)=d(x,y)\quad\textrm{and}\quad d(\delta_{r}(x),\delta_{r}(y))=r\,d(x,y)\quad\textrm{for all }x,y,z\in\mathbb G,\, r>0.
\]
We denote by $\mathcal H^q$ and $\mathcal S^q$ the $q$-dimensional Hausdorff and spherical Hausdorff 
measures induced by a fixed homogeneous distance $\rho$, respectively. Open balls with center at $x$ 
and radius $r>0$ with respect to $d$ will be denoted by $B_{x,r}$ and the corresponding closed balls
will be denoted by $D_{x,r}$.

\subsection{Graded groups as vector spaces}\label{subgv}

The fact that the exponential mapping $\exp:\cG\lra\G$ is a diffeomorphism
allows us to identify $\G$ with its Lie algebra as follows. 
We set $H_j=\exp V_j$ for all $j=1,\ldots,\iota$ and equip $\G$ with a structure of 
graded vector space satisfying $\G=H_1\oplus\cdots \oplus H_\iota$.
This direct sum has the corresponding {\em canonical projections}
\[
p_j:\G\lra H_j,\quad p_j(\sum_{l=1}^\iota x_l)=x_j,\quad \mbox{where}\quad
x_l\in H_l\quad \mbox{for all $l=1,\ldots,\iota$}\,.
\]
The group operation 
\begin{equation}\label{groupope}
x\cdot y=x+y+Q(x,y)
\end{equation}
has a polynomial form, where $Q$ is given by the Baker-Campbell-Hausdorff formula. 
We define the integers $m_0=0$ and $m_j=\sum_{l=1}^j\dim H_l$ for any
$j=1,\ldots,\iota$. Then a {\em graded basis} $(e_1,\ldots,e_n)$ of $\G$ is one
of its bases such that 
\[
(e_{m_{j-1}+1},\ldots,e_{m_j})\quad\mbox{is a basis of}\quad H_j
\]
for all $j=1,\ldots,\iota$. In the sequel, a graded basis will be understood 
whenever $\G$ is identified with $\R^n$. Declaring this basis orthonormal, we have fixed 
both a scalar product on $\G$ and a left invariant metric $g$ on $\G$ with respect to
the group operation.
With slight abuse of notation, we use the same notation to denote both the length $|v|$ 
of a vector $v\in T_y\G$ with respect to $g$ and to denote the 
norm $|x|$ of an element $x$ of $\G$. Notice that the metric $g$ at the origin 
exactly coincides with the fixed scalar product of $\G$, since $g$ is defined as 
the unique left invariant Riemannian metric on $\G$ with this property.
All the graded bases we consider are understood to be orthonormal with respect to the
underlying scalar product $g$.

The basis $(e_l)$ along with its coordinates $(x_l)$ automatically inherit a {\em degree}
from the layers $H_j$, namely, 
\[
d_j=k \quad\mbox{if and only if} \quad e_j\in H_k\,.
\]
Elements of $H_k$ have degree $k$ and any coordinate $x_j$ of an element
$\sum_{j=1}^nx_je_j\in\G$ has degree $d_j$.
Taking into account this notion, dilations on $\G$ can be written as
\begin{equation}\label{dilat}
\delta_r\Big(\sum_{l=1}^nx_l\,e_l\Big)=\sum_{l=1}^n r^{d_l}\,x_l\,e_l\,.
\end{equation}

Notice that a change of system of graded coordinates is a linear isomorphism,
hence the class of polynomials on $\G$ is well defined. In fact, a polynomial on $\G$ 
is simply a standard polynomial on $\R^n$ when some fixed graded coordinates are fixed on $\G$.
Let us stress that a linear isomorphism $I:\R^n\lra\R^n$ arising from a change 
of graded coordinates automatically preserves the degrees of all coordinates, namely, 
\[
I(e_j)\in\span\{e_l: l=m_{d_j-1}+1, \ldots, m_{d_j}\} \quad \mbox{ for all}
\quad m_{d_j-1}+1\leq j\leq m_{d_j}.
\]
For any $\alpha\in\N^n$ we define the monomial 
\[
x^\alpha:=x_1^{\alpha_{1}}x_2^{\alpha_{2}}\ldots x_n^{\alpha_n}\quad\textrm{where} \quad 
d(\alpha):=\sum_{i=1}^n d_j\,\alpha_j
\]
is the {\em homogeneous degree} of $x^\alpha$. By \eqref{dilat} we have
\[
(\delta_rx)^\alpha=r^{d(\alpha)}\prod_{j=1}^nx_j,
\]
where $\alpha=(\alpha_1,\ldots,\alpha_n)$. 
A polynomial is {\em $d$-homogeneous} if it is a linear combination of monomials 
of homogeneous degree equal to $d$. Thus, a polynomial $P$ is $d$-homogeneous if and 
only if $P(\delta_rx)=r^d\,P(x)$ for all $x\in\R^n$ and $r>0$. In this case, 
we also say that the {\em homogeneous degree} of $P$ is $d$. 

The degree transmits itself to the unique left invariant vector fiels $X_j$ of $\G$
such that $X_j(0)=e_j$. Then $X_j$ has degree $d_j$ and has the polynomial form
\begin{equation}\label{eq2}
X_j(x)=\der_{x_j}+\sum_{l:\,d_l>d_j} a^l_j(x)\,\der_{x_l}
\end{equation}
with respect to the understood graded coordinates $(x_j)$.
The polynomails $a^l_j$ in the previous formula are $(d_l-d_j)$-homogeneous.
If $(X_1,\ldots,X_n)$ is a frame of left invariant vector fields 
generated by a graded basis $(e_1,\ldots,e_n)$ of $\G$, then we will say that
this frame is also {\em graded}.

Next, we introduce the notion of degree of curves in a graded group.
We address the reader to the original work, \cite{MagVit}, for more details and the
more general notion of degree of a submanifold. 
\begin{Def}[Degree of curves]\rm
Let $\Sigma$ be a 1-dimensional submanifold of a graded group $\G$ and let $x\in\Sigma$. 
Let $\tau\in T_x\Sigma\sm\{0\}$ and let $\tau=\sum_{j=1}^n\mu_j\,X_j(x)$.
The {\em pointwise degree} of $\Sigma$ at $x$ is the integer
$d_\Sigma(x)=\max\{d_j: \lambda_j\neq0 \}.$ The degree of $\Sigma$ is given by
\[d(\Sigma)=\max_{x\in\Sigma}d_\Sigma(x).\]
\end{Def} 
The left translations of the group are denoted by $l_x:\G\lra\G$, $l_x(y)=x\cdot y$.
\begin{Def}\rm 
Let $\tilde g$ be a Riemannian metric in $\G$, let $\Sigma$ be a one dimensional submanifold of $\G$
and let $\gamma:(a,b)\lra\G$ be a $C^1$ smooth local parametrization of $\Sigma$, with $t\in(a,b)$ 
and $x=\gamma(t)$.
Then we define the {\em unit tangent vector} of $\Sigma$ at $x$ as 
$\tau_\Sigma(x)=\dot\gamma(t)/|\dot\gamma(t)|_{\tilde g}\in T_x\Sigma$.
If $p_j:\G\lra H_j$, $1\leq j\leq\iota$ is the canonical projection of $\G$, 
and we have identified $\G$ with $T_0\G$, then we define the
{\em $j$-projection of $\tau_\Sigma(x)$} by
\begin{equation}\label{protau}
\tau_\Sigma^j(x)=dl_{x}\big(p_j(dl_{x^{-1}}\big(\tau_\Sigma(x)\big)\big)\in T_x\G\cap dl_x(H_j)\,.
\end{equation}
\end{Def}
\begin{Def}[Metric factor]\label{metricf}\rm
Let $x\in\G$, let $\tau\in T_x\G$ and set $\tau_0=dl_{x^{-1}}\tau\in T_0\G$.
Then the {\em metric factor} is defined by
\begin{eqnarray}\label{metricfactor}
\theta(\tau)=\cH_{|\cdot|}^1\left(\span\{\tau_0\}\cap B_1\right),
\end{eqnarray}
where $|\cdot|$ is the fixed scalar product on $\G$ and $B_1$ is the open unit ball with respect
to the fixed homogeneous distance $d$.
\end{Def}
Clearly, the metric factor is constant on directions of a left invariant vector field.

\subsection{Some density estimates of Geometric Measure Theory} 

Our arguments are based on the following elementary fact, 
see for instance 2.10.19 of \cite{Fed}.
\begin{Lem}
Let $X$ be a metric space, let $\mu$ be a Borel measure on $X$
and let $\{V_i\}_{i\in\N}$ be an open covering of $X$
such that $\mu(V_i)<\infty$. 
Let $Z\subset X$ be a Borel set and
suppose that
\[
\limsup_{r\to0^+}r^{-a}\mu(D_{x,r})\geq\kappa>0
\]
whenever $x\in Z$, where $a>0$. Then $\mu(Z)\geq\kappa\;\cS^a(Z)$.
\end{Lem}
The symbol $\cS^a$ in the previous lemma denotes the $a$-dimensional spherical Hausdorff measure
constructed by the size function $\zeta_a(D_{x,r})=r^a$ and
$D_{x,r}$ is the closed ball of center $x$ and radius $r$.
Then we have the following 
\begin{Cor}\label{gmta}
Let $\Sigma$ be $1$-dimensional $C^{1}$ submanifold of a stratified group $\mathbb G$ and
let $\mu$ be the left invariant Riemannian measure of $\G$
restricted to $\Sigma$. If $Z$ is a Borel set of $\Sigma$ such that
$\limsup_{r\to0^+}r^{-a}\mu(D_{z,r})=+\infty$, whenever
$z\in Z$, then $\cS^a(Z)=0$.
\end{Cor}
\section{Blow-ups and negligibility}

This section is devoted to the proofs of our main results.

\begin{Pro}\label{blw}
Let $\gamma:(-1,1)\lra\G$ be $C^1$ smooth embedding, let $\Sigma$ be its image
and fix $d(\Sigma)=q$. If $\gamma(0)=0$ and $d_\Sigma(0)=q$, then there exists a 
graded basis $(e_1,\ldots,e_n)$ such that for all $j=1,2,\ldots,n$ we have
\begin{equation}\label{estblw}
\gamma_i(t)=o(t^{d_i/q}) \quad\mbox{whenever}\quad i\neq i_0=m_{q-1}+1.
\end{equation}
\end{Pro}
\begin{proof}
The grading $\G=H_1\oplus \cdots\oplus H_\iota$ has the associated
canonical projections $p_j:\G\lra H_j$, hence $p_q(\dot\gamma(0))\neq0$.
We set $e_{i_0}=p_q(\dot\gamma(0))/|p_q(\dot\gamma(0))|$.
We complete this vector to a graded basis $(e_1,\ldots,e_n)$ of $\G$
and consider the corresponding frame of left invariant vector fields
$(X_1,\ldots,X_n)$. We have 
\[
\dot\gamma(t)=\sum_{s=1}^n C^s(\gamma(t))\,X_s(\gamma(t))\,,
\]
where the functions $t\to C^s(\gamma(t))$ are continuous, vanish on a neighbourhood of $0$
whenever $d_s>q$, $C^{i_0}(0)\neq0$ and for all $s$ such that
$i_0<s\leq m_q$ we have $C^s(0)=0$. 
Let us introduce the homeomorphism $\eta:\R\lra\R$ defined as 
$\eta(t)=\big(|t|^q\,\mbox{sgn}(t)\big)/q$.
We consider the reparametrized curve $\sigma(t)=\gamma\big(\eta(t)\big)$, hence
\[
\sigma'(t)=|t|^{q-1} \sum_{s:\,d_s\leq q} C^s(\sigma(t))\,X_s(\sigma(t))\,.
\]
Thus, clearly $\sigma_i'(t)=o(t^{d_i-1})$ for all $i$ such that $d_i<q$. 
Due to \eqref{eq2}, setting $a^i_s\equiv\delta^i_s$ when $d_i=d_s$ 
and $a^i_s\equiv 0$ when $d_i<d_s$, we can write
\begin{equation*}
X_j(x)=\der_{x_j}+\sum_{l=1}^n a^l_j(x)\,\der_{x_l}\,.
\end{equation*}
We consider the cases $d_i\geq q$, where the previous formula gives
\begin{equation}\label{q-1}
\sigma'_i(t)=|t|^{q-1}\sum_{s:\,d_s\leq q} C^s(\sigma(t))\, a^i_s(\sigma(t))\,.
\end{equation}
We first consider all $i$'s such that $d_i=q$, therefore  
\begin{eqnarray*}
\sigma'_i(t)&=&|t|^{q-1}\sum_{s:\,d_s\leq q} C^s(\sigma(t))\, a^i_s(\sigma(t))\\
&=&|t|^{q-1}\,C^i(\sigma(t))+|t|^{q-1}\sum_{s:\,d_s<q}C^s(\sigma(t))\,a^i_s(\sigma(t))\,.
\end{eqnarray*}
The polynomials $a^i_s(x)$ are $(q-d_s)$-homogeneous, hence they only depend on 
the components $x_l$ with $d_l<q$, for which $\gamma_l(t)=o(t^{d_l})$. 
Then $a^i_s(\sigma(t))=o(t^{q-d_s})$ and
\begin{eqnarray*}
\sigma'_i(t)=|t|^{q-1}\,C^i(\sigma(t))+o(t^{q-1})\,.
\end{eqnarray*}
In the case $i\neq i_0$ and $d_i=q$, we have $C^i(0)=0$, therefore $\sigma'_i(t)=o(t^{q-1})$.
This is not true for $i=i_0$, since $C^{i_0}(0)\neq0$.

Now, we cosider all $i$'s such that $d_i=q+1$. We split \eqref{q-1} into two addends, getting
\[
\sigma_i'(t)=|t|^{q-1}\sum_{s:\,d_s=q} C^s(\sigma(t))\,a^i_s(\sigma(t))
+|t|^{q-1}\sum_{s:\,d_s<q} C^s(\sigma(t))\,a^i_s(\sigma(t))\,.
\]
In the second addend, the homogeneous degree of $a^i_s$ is less than or equal to $q$,
hence we can only conclude that $a^i_s(\gamma(t))=O(t^{q+1-d_s})=o(t)$. It follows that
\[
\sigma_i'(t)=|t|^{q-1}\sum_{s:\,d_s=q} C^s(\sigma(t))\,a^i_s(\sigma(t))+o(t^q).
\]
Furthermore, $C^s(\sigma(t))=o(1)$ whenever $s\neq i_0$ and $d_s=q$
and $a^i_s(\sigma(t))=O(t)$. We have then established
\begin{equation}\label{lastal}
\sigma_i'(t)=|t|^{q-1}\,C^{i_0}(\sigma(t))\,a^i_{i_0}(\sigma(t))+o(t^q).
\end{equation}
Now, by definition of left invariant vector fields, and representing 
$Q(x,y)$ of \eqref{groupope} as $\sum_{i=1}^n Q_i(x,y)\,e_i$, we have 
$a^i_s(x)=\der_{y_s}Q_i(x,0)$ whenever $d_i>d_s$.
Clearly $\span\{X_{i_0}\}$ is a subalgebra of $\cG$, then we can apply
Lemma~2.5 of \cite{MagVit}, that gives
\[
Q_i(x,y)=\sum_{l\neq i_0,\,d_l<d_i} x_l\, R^i_l(x,y)+y_l\,S^i_l(x,y)
\]
whenever $i\neq i_0$. In particular, we are lead to
\begin{equation}\label{algeq}
a^i_{i_0}(x)=\sum_{l\neq i_0,\,d_l<d_i} x_l\; \der_{i_0}R^i_l(x,0)\,.
\end{equation}
In the previous steps, we have proved that $\sigma_l(t)=o(t^{d_l})$ for each $l$ such that 
$d_l\leq q$ and $l\neq i_0$, hence the previous formula yields
\begin{equation*}
a^i_{i_0}(\sigma(t))=\sum_{l\neq i_0,\,d_l<d_i} \sigma_l(t)\; 
\der_{i_0}R^i_l(\sigma(t),0)=o(t^{d_i-d_{i_0}})\,,
\end{equation*}
where we have taken into account that $\der_{i_0}R^i_l(\sigma(t),0)=O(t^{d_i-d_0-d_l})$
for $d_i-d_{i_0}-d_l\geq0$ and $\der_{i_0}R^i_l(\sigma(t),0)\equiv 0$ otherwise.
In view of \eqref{lastal}, we have proved that $\sigma_i'(t)=o(t^q)$.

To complete the proof, we argue by induction, assuming that 
$\sigma_i'(t)=o(t^{d_i-1})$ whenever $i$ satisfies $1+q\leq d_i\leq k-1$ and $k\geq 2$.
Then we consider the cases $d_i=k$. 
The next steps essentially repeat the previous argument. In fact, we now consider
\[
\sigma_i'(t)=|t|^{q-1}\sum_{s:\,d_s=q} C^s(\sigma(t))\,a^i_s(\sigma(t))
+|t|^{q-1}\sum_{s:\,d_s<q} C^s(\sigma(t))\,a^i_s(\sigma(t))\,.
\]
The homogeneous degree of $a^i_s$ is less than or equal to $k-1$,
hence $a^i_s(\gamma(t))=O(t^{d_i-d_s})$. It follows that
\[
\sigma_i'(t)=|t|^{q-1}\sum_{s:\,d_s=q} C^s(\sigma(t))\,a^i_s(\sigma(t))+o(t^{d_i-1}).
\]
Since $C^s(0)=0$ for $d_s=q$ and $s\neq i_0$, we have 
\[
C^s(\sigma(t))\,a^i_s(\sigma(t))=o(1)\,O(t^{d_i-d_s})=o(t^{d_i-d_s}).
\]
It follows that 
$\displaystyle
\sigma_i'(t)=|t|^{q-1}\, C^{i_0}(\sigma(t))\,a^i_{i_0}(\sigma(t))
+o(t^{d_i-1})\,.
$
Using formula \eqref{algeq} and the induction hypothesis, we also get
$a^i_{i_0}(\sigma(t))=o(t^{d_i-d_{i_0}})$, leading us to $\sigma'_i(t)=o(t^{d_i-1})$,
that concludes the proof.
\end{proof}
The previous proposition is the basic ingredient to establish the blow-up
at points of maximum degree.
\begin{proof}[Proof of Theorem~\ref{blwmea}]
For any sufficiently small $r>0$, we have
\begin{eqnarray}
\frac{\tilde\mu_1(\Sigma\cap B_{x,r})}{r^q} &=& 
\frac{1}{r^q}\int_{\Gamma^{-1}(B_{x,r})} |\dot\Gamma(t)|_{\tilde g}\,dt \nonumber\\
&=& \int_{A_{x,r}} |\dot\Gamma(r^qt)|_{\tilde g}\; dt\,,
\end{eqnarray}
where $\Gamma:(-1,1)\lra\G$ is a local chart of $\Sigma$,
$\Gamma(0)=x$ and $\Sigma\cap B_{x,r}=\Gamma\big(\Gamma^{-1}(B_{x,r})\big)$.
We have also set $A_{x,r}=\{t\in(-1,1):  \Gamma(tr^q)\in B_{x,r}\}.$
Let us consider the curve $\gamma(t)=x^{-1}\cdot\Gamma(t)$ and observe that
it parametrizes the translated submanifold $x^{-1}\Sigma$, that has same degree $q$,
and $\gamma(0)=0$. Then 
\[
A_{x,r}=\left\{t\in(-r^{-q},r^{-q}): \sum_{i=1}^n\frac{\gamma_i(tr^q)}{r^{d_i}}\,e_i\in B_1\right\}
\]
and Proposition~\ref{blw} yields a graded basis $(e_1,\ldots,e_n)$ with respect to
which we have $\gamma_i(t)=o(|t|^{d_i/q})$ whenever $i\neq i_0$ and $i_0=m_{q-1}+1$.
Then for these integers
\[
\frac{\gamma_i(tr^q)}{r^{d_i}}=\frac{o(|tr^q|^{d_i/q})}{r^q}\lra 0\quad\mbox{as}\quad r\to0^+
\]
and clearly $\gamma_{i_0}(tr^q)r^{-q}\to\gamma_{i_0}'(0)\,t$, where these limits are
uniform with respect to $t$ that varies on compact sets of $\R$.
In particular, we have the $L^1_{loc}$ convergence of the characteristic functions
\[
{\bf 1}_{A_{x,r}}\lra{\bf 1}_{\frac{1}{\gamma_{i_0}'(0)}S_0}, \quad\mbox{where}\quad 
\displaystyle S_0=\big\{t\in\R: te_{i_0}\in B_1\big\}
\]
and ${\bf 1}_{A_{x,r}}\leq{\bf 1}_{[-M_0,M_0]}$ for some $M_0>0$ and any $r>0$  small.
By Lebesgue's Theorem,
\begin{equation}\label{blwlim}
\lim_{r\to0^+}\frac{\tilde\mu_1(\Sigma\cap B_{x,r})}{r^q}=|\dot\Gamma(0)|
\cL^1\left(\frac{1}{\dot\gamma_{i_0}(0)}S_0\right)=
\frac{|\dot\Gamma(0)|}{|\dot\gamma_{i_0}(0)|}\,\cL^1(S_0).
\end{equation}
We first observe that $|e_{i_0}|=1$ gives $\cL^1(S_0)=\cH^1_{|\cdot|}\big(B_1\cap\span\{e_{i_0}\}\big)=
\theta(\tau_\Sigma^q(x))$.
Then we observe that 
\[
 \frac{\dot\gamma_{i_0}\,e_{i_0}}{|\dot\Gamma(0)|_{\tilde g}}
=\frac{p_q\big(dl_{x^{-1}}\dot\Gamma(0)\big)}{|\dot\Gamma(0)|_{\tilde g}}
=dl_{x^{-1}}\left(\frac{dl_x\big(p_q\big(dl_{x^{-1}}\dot\Gamma(0)\big)\big)}{|\dot\Gamma(0)|_{\tilde g}}\right)
=dl_{x^{-1}}\big(\tau^q_\Sigma(x)\big)\,,
\]
hence we have 
\[
|\tau^q_\Sigma(x)|=|dl_{x^{-1}}\big(\tau^q_\Sigma(x)\big)|=
\left|\frac{\dot\gamma_{i_0}\,e_{i_0}}{|\dot\Gamma(0)|_{\tilde g}}\right|=
\frac{|\dot\gamma_{i_0}|}{|\dot\Gamma(0)|_{\tilde g}}.
\]
The last equalities joined with \eqref{blwlim} lead us to \eqref{1blw}, concluding the proof.
\end{proof}
\begin{Pro}\label{neglestim}
Let $\gamma:(-1,1)\lra\G$ be $C^1$ smooth embedding and
let $\Sigma$ be the image of $\gamma$, namely, a 1-dimensional submanifold of $\G$.
Assuming that $\gamma(0)=0$ and $d_\Sigma(0)<d(\Sigma)$, then for every
graded basis $(e_1,\ldots,e_n)$ with $\gamma=\sum_{j=1}^n\gamma_j\,e_j$, we have nonnegative 
functions $\ep_j$ such that for all $t$ small and any $j=1,2,\ldots,n$, we have
\begin{equation}\label{claim}
|\gamma_j(t)|\leq \varepsilon_j(|t|)\,|t|^{p_j},
\end{equation}
where $p_j=d_j/d(\Sigma)$ and $\varepsilon_j(r)\to 0^+$ as $r\to0^+$.
\end{Pro}
\begin{proof}
There exists a neighbourhood $I\subset \mathbb R$ of $0$ such that 
\begin{equation}\label{eq1}
\dot{\gamma}(t)=\sum_{j=1}^n \lambda_j(t)X_j(\gamma(t))
=\sum_{j=1}^n \dot{\gamma}_j(t)\,\partial_{x_j}.
\end{equation}
We observe that the second equality and the $C^1$ smoothness of $\gamma$ imply that 
all $\lambda_j$ are continuous.
Let us first consider the case $d_j \leq d_\Sigma(0)$.  Then $p_j<1$ and  
$C^1$ continuity of $\gamma_j$ yields  
$| \dot {\gamma_{j}}(t) |< \max\{2 | \dot{ \gamma_{j}}(0) |,1\} $ for $|t|$ small enough.
Therefore, as $\gamma_j(0)=0$, we have $| \gamma_j (t)|< \max\{2 |\dot{\gamma_j}(0)|,1\}|t|$.
This proves the validity of \eqref{claim}, where one sets
$\varepsilon_{j}(r)=\max\{2|\dot{\gamma_j}(0)|,1\} r^{1-p_{j}}$.
If $d_\Sigma(0)<d_j\leq d(\Sigma)$, then we must have $\dot{\gamma_j}(0)=0$. 
Thus, introducing the nondecreasing function $f_j(r):=\max_{|s|\leq r}|\dot{\gamma_{j}}(s)|$
we have
\[
\lim_{r\to0^+} f_j(r)=0\quad\textrm{and}\quad |\dot{\gamma_{j}}(t)|\leq f_{j}(t).
\]
Therefore $|\gamma_{j}(t)|\leq f_{j}(|t|) |t| = \varepsilon_{j}(|t|)|t|^{p_{j}}$,
where we have set $\varepsilon_{j}(r)=f_{j}(r)\,r^{1-p_{j}}$. 
Now, to prove the validity of \eqref{claim} for all $j=1,\ldots,n$, we use
the following induction argument.
Suppose that the for every $n$-dimensional stratified group $\G$ and every $\Sigma$ and $\gamma$ 
satisfying our assumptions, the claim \eqref{claim} holds for all $j=1,2,\ldots,k-1$, where $k\leq n$. 
Then we have to show that this yields \eqref{claim} for $j=k$.

If $k\in\N$ is such that $d_k\leq d(\Sigma)$, then \eqref{claim} holds,
due to the previous arguments. Let us consider the case $d_k>d(\Sigma)$.
Joining \eqref{eq2} with \eqref{eq1}, we have 
\[
\dot{\gamma_k}(t)=\lambda_k(t)+\sum_{d_j<d_k}\lambda_j(t)a^k_j(\gamma(t)).
\]
By definition of degree, if $d_j>d(\Sigma)$, then 
$\lambda_j$ everywhere vanishes. Therefore it remains to estimate $a_j^k\circ \gamma$ 
whenever $d_j\leq d(\Sigma)$.
It follows that  
\[
\dot{\gamma_k}(t)=\sum_{j:\,d_j\leq d(\Sigma)}\lambda_j(t)a^k_j(\gamma(t)).
\]
Recall that $a^k_j$ is a polynomial $(d_k-d_j)$-homogeneous. Then it only depends on variables
$x_l$'s with $d_l<d_k$. Since $l\to d_l$ is nondecreasing, $a^k_j$ only depends on the
variables $x_l$ with $l<k$. It follows that $a^k_j$ can be written as 
\[
\sum_{d(\alpha)=d_k-d_j}c_\alpha\,x^\alpha=\sum_{d(\alpha)=d_k-d_j} c_\alpha\;x_1^{\alpha_1}\cdots x_{k-1}^{\alpha_{k-1}}\,.
\]
Taking into account the induction hypothesis, we get 
\[
\Big|\sum_{d(\alpha)=d_k-d_j}c_\alpha\,\gamma(t)^\alpha\Big|\leq
\bigg(\sum_{d(\alpha)=d_k-d_j}c_\alpha\;\prod_{l=1}^{k-1}\ep_j(t)^{\alpha_l}\bigg)\; |t|^{\sum_{l=1}^{k-1}p_l\alpha_l}.
\]
We also observe that 
\[
\sum_{l=1}^{k-1}p_l\,\alpha_l=\sum_{l=1}^{k-1}\frac{d_l}{d(\Sigma)}\,\alpha_l=\frac{d_k-d_j}{d(\Sigma)}
\geq p_k-1\,.
\]
It follows that 
\[
\Big|\sum_{d(\alpha)=d_k-d_j}c_\alpha\,\gamma(t)^\alpha\Big|\leq
\bigg(\sum_{d(\alpha)=d_k-d_j}c_\alpha\;\prod_{l=1}^{k-1}\ep_j(|t|)^{\alpha_l}\bigg)\; |t|^{p_k-1}.
\]
Then we can find a nonnegative and nondecreasing function $\ep_k(r)$ for $r\geq0$ sufficiently
small,  only depending on all polynomials $a^k_j$
and all $\ep_j$ for $j\leq k-1$ such that 
\[
|\dot\gamma_k(t)|\leq\ep_k(|t|)\, |t|^{p_k-1}\quad\mbox{and}\quad \lim_{r\to0^+}\ep_k(r)=0\,.
\] 
This immediately leads us to the validity of \eqref{claim} for $j=k$ and completes the proof.
\end{proof}
\begin{Pro}\label{limitinfinite} Let $\Sigma$ be a $C^1$ smooth $1$-dimensional submanifold in a stratified group $\mathbb G$. 
Let $\mu$ be the Riemannian measure induced on $\Sigma$ by to the fixed left invariant metric $g$.
Then for every $x\in \Sigma$ such that $d_{\Sigma}(x)<d(\Sigma)$, there holds
\[
\lim_{r\rightarrow 0^+}\frac{\mu(\Sigma\cap D_{x,r})}{r^{d(\Sigma)}}=+\infty.
\]
\end{Pro}
\begin{proof}
For each $j=1,\ldots,n$, we define the orthogonal projection $\mathfrak p_j:\G\lra\span\{e_j\}$
associated to the graded basis $(e_1,\ldots,e_n)$ of $\G$. We define the box
\[
\Bx_r=\{x\in\G: |\cp_j(x)|\leq r^{d_j}\}.
\]
Since $\delta_r(\Bx_1)=\Bx_r$ for all $r>0$ and $\Bx_\lambda\subset D_1\subset \Bx_{\lambda^{-1}}$
for a fixed $\lambda>0$, we get 
\begin{equation}\label{ballbox}
\Bx_{r\lambda}\subset D_r\subset \Bx_{r/\lambda}\,.
\end{equation}
Let $\gamma:(-1,1)\lra\G$ be a $C^1$ smooth embedding such that
$\gamma(0)=0$ and whose image is contained in the translated
submanifold $x^{-1}\cdot\Sigma$. Let $q=d(\Sigma)$ and consider
\[
\frac{\mu(\Sigma\cap D_{x,r})}{r^q}=\frac{\mu\big((x^{-1}\cdot\Sigma)\cap D_r\big)}{r^q}
\geq \frac{\mu\big((x^{-1}\cdot\Sigma)\cap\Bx_{\lambda r}\big)}{r^q}\,.
\]
Then for $r>0$ sufficiently small, we have
\begin{equation}
\frac{\mu(\Sigma\cap D_{x,r})}{r^q}\geq 
\frac{1}{r^q}\int_{I_{x,r}}|\dot\gamma(t)|\,dt.
\end{equation}
In the previous inequality we have set $I_{x,r}=\displaystyle \gamma^{-1}(\Bx_{\lambda r})=\bigcap_{j=1}^n A_j(x,r)$, where
\[ 
A_j(x,r)=\{t\in(-1,1): |\gamma_j(t)|\leq (\lambda r)^{d_j}\}\,.
\]
Let $0<h<1$ be arbitrarily fixed and let $0<\delta_j<1$ be such that 
$\ep_j(|t|)<h$ whenever $|t|\leq\delta_j$.
Let $\ep_j$ be as in Proposition~\ref{neglestim} and fix an arbitrary $h>0$.
Then there exist $\delta_j>0$ such that $\ep_j(|t|)\leq h$ whenever $|t|\leq\delta_j$.
It follows that 
\[
\{t\in(-\delta_j,\delta_j): h\,|t|^{d_j/q}\leq (\lambda r)^{d_j}\}\subset A_j(x,r)
\]
for all $j=1,\ldots,n$. Now, we define the numbers
\[
0<\sigma_0<\min\{\delta_j: j=1,\ldots,n\}\quad\mbox{and}\quad
\sigma_1=\max\{\delta_j: j=1,\ldots,n\}
\]
such that $|\dot\gamma(t)|>|\dot\gamma(0)|/2$ for any $t\in(-\sigma_0,\sigma_0)$. Then we introduce the set  
\[
S_r=\{t\in(-\sigma_0,\sigma_0): |t| \leq (\lambda^q\,r^q)/h^{q/\sigma_1}\}.
\]
We observe that $S_r\subset A_j(x,r)$ for every $j=1,\ldots,n$ and for $r>0$
sufficiently small
\[
\frac{\mu(\Sigma\cap D_{x,r})}{r^q}\geq 
\frac{|\dot\gamma(0)|}{2r^q}|S_r|=\frac{|\dot\gamma(0)|}{2}\frac{\lambda^q}{h^{q/\sigma_1}}\,.
\]
Taking into account the arbitrary choice of $h$, we have reached our claim. 
\end{proof}

\begin{proof}[Proof of Theorem~\ref{negconcur}]
Taking into account Corollary~\ref{gmta}, the proof is an immediate consequence
of Proposition~\ref{limitinfinite}.
\end{proof}

\end{document}